# PCA CONSISTENCY IN HIGH DIMENSION, LOW SAMPLE SIZE CONTEXT

By Sungkyu Jung and J. S. Marron

*University of North Carolina*

Principal Component Analysis (PCA) is an important tool of dimension reduction especially when the dimension (or the number of variables) is very high. Asymptotic studies where the sample size is fixed, and the dimension grows [i.e., High Dimension, Low Sample Size (HDLSS)] are becoming increasingly relevant. We investigate the asymptotic behavior of the Principal Component (PC) directions. HDLSS asymptotics are used to study consistency, strong inconsistency and subspace consistency. We show that if the first few eigenvalues of a population covariance matrix are large enough compared to the others, then the corresponding estimated PC directions are consistent or converge to the appropriate subspace (subspace consistency) and most other PC directions are strongly inconsistent. Broad sets of sufficient conditions for each of these cases are specified and the main theorem gives a catalogue of possible combinations. In preparation for these results, we show that the geometric representation of HDLSS data holds under general conditions, which includes a $\rho$-mixing condition and a broad range of sphericity measures of the covariance matrix.

**1. Introduction and summary.** The High Dimension, Low Sample Size (HDLSS) data situation occurs in many areas of modern science and the asymptotic studies of this type of data are becoming increasingly relevant. We will focus on the case that the dimension $d$ increases while the sample size $n$ is fixed as done in Hall, Marron and Neeman [8] and Ahn et al. [1]. The $d$-dimensional covariance matrix is challenging to analyze, in general, since the number of parameters is $\frac{d(d+1)}{2}$, which increases even faster than $d$. Instead of assessing all of the parameter estimates, the covariance matrix is usually analyzed by Principal Component Analysis (PCA). PCA is often used to visualize important structure in the data, as shown in Figure 1. The data in Figure 1, described in detail in Bhattacharjee et al. [4] and Liu et al. [15], are from a microarray study of lung cancer. Different symbols correspond to cancer subtypes, and Figure 1 shows the projections of the data onto the subspaces generated by PC1 and PC2 (left panel) and PC1 and PC3 (center panel, resp.) directions. This shows the difference between subtypes is so strong that it drives the first three principal components. This illustrates a common occurrence: the data have an important underlying structure which is revealed by the first few PC directions.

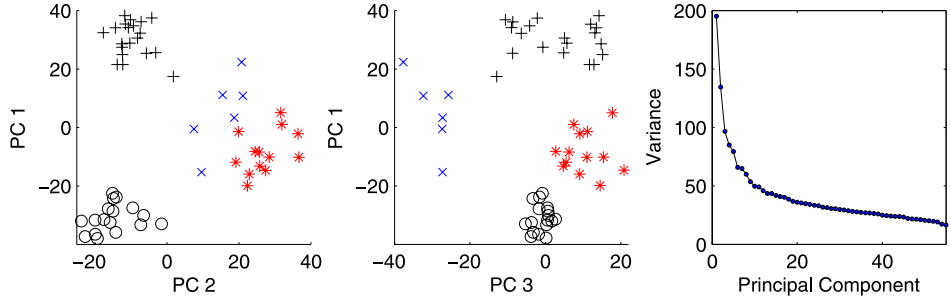

FIG. 1. *Scatterplots of data projected on the first three PC directions. The dataset contains 56 patients with 2530 genes. There are 20 Pulmonary Carcinoid (plotted as +), 13 Colon Cancer Metastases (∗), 17 Normal Lung (◦), and 6 Small Cell Carcinoma (×). In spite of the high dimensionality, PCA reveals important structure in the data. This corresponds to the* consistent *case in our asymptotics, as shown in the scree plot on the right. Note that the first few eigenvalues are much larger than the rest.*

PCA is also used to reduce dimensionality by approximating the data with the first few principal components.

For both visualization and data reduction, it is critical that the PCA empirical eigenvectors reflect true underlying distributional structure. Hence, our focus is on the underlying mechanism which determines when the sample PC directions converge to their population counterparts as $d \to \infty$. In general, we assume $d > n$. Since the size of the covariance matrix depends on $d$, the population covariance matrix is denoted as $\Sigma_d$ and similarly the sample covariance matrix, $S_d$, so that their dependency on the dimension is emphasized. PCA is done by eigen decomposition of a covariance matrix. The eigen decomposition of $\Sigma_d$ is

$$\Sigma_d = U_d \Lambda_d U_d',$$

where $\Lambda_d$ is a diagonal matrix of eigenvalues $\lambda_{1,d} \geq \lambda_{2,d} \geq \cdots \geq \lambda_{d,d}$ and $U_d$ is a matrix of corresponding eigenvectors so that $U_d = [u_{1,d}, u_{2,d}, \ldots, u_{d,d}]$. $S_d$ is similarly decomposed as

$$S_d = \hat{U}_d \hat{\Lambda}_d \hat{U}_d'.$$

Ahn et al. [1] developed the concept of HDLSS consistency which was the first investigation of when PCA could be expected to find important structure in HDLSS data. Our main results are formulated in terms of three related concepts:

1. *Consistency*: The direction $\hat{u}_{i,d}$ is *consistent* with its population counterpart $u_{i,d}$ if $\text{Angle}(u_{i,d}, \hat{u}_{i,d}) \to 0$ as $d \to \infty$. The growth of dimension can be understood as adding more variation. The consistency of sample eigenvectors occurs when the added variation supports the existing structure in the covariance or is small enough to be ignored.

2. *Strong inconsistency*: In situations where $\hat{u}_{i,d}$ is not consistent, a perhaps counter-intuitive HDLSS phenomenon frequently occurs. In particular, $\hat{u}_{i,d}$ is said to be *strongly inconsistent* with its population counterpart $u_{i,d}$ in the sense that it tends to be as far away from $u_{i,d}$ as possible, that is, $\text{Angle}(u_{i,d}, \hat{u}_{i,d}) \to \frac{\pi}{2}$ as $d \to \infty$. Strong inconsistency occurs when the added variation obscures the underlying structure of the population covariance matrix.

3. *Subspace consistency*: When several population eigenvalues indexed by $j \in J$ are similar, the corresponding sample eigenvectors may not be distinguishable. In this case, $\hat{u}_{j,d}$ will not be consistent for $u_{j,d}$ but will tend to lie in the linear span, $\text{span}\{u_{j,d} : j \in J\}$. This motivates the definition of convergence of a direction $\hat{u}_{i,d}$ to a subspace, called *subspace consistency*;

$$\text{Angle}(\hat{u}_{i,d}, \text{span}\{u_{j,d} : j \in J\}) \longrightarrow 0$$

as $d \to \infty$. This definition essentially comes from the theory of *canonical angles* discussed by Gaydos [7]. That theory also gives a notion of convergence of subspaces, that could be developed here.

In recent years, substantial work has been done on the asymptotic behavior of eigenvalues of the sample covariance matrix in the limit as $d \to \infty$, see Baik, Ben Arous and Péché [2], Johnstone [11] and Paul [16] for Gaussian assumptions and Baik and Silverstein [3] for non-Gaussian results when $d$ and $n$ increase at the same rate, that is, $\frac{n}{d} \to c > 0$. Many of these focus on the *spiked covariance model*, introduced by Johnstone [11]. The spiked covariance model assumes that the first few eigenvalues of the population covariance matrix are greater than 1 and the rest are set to be 1 for all $d$. HDLSS asymptotics, where only $d \to \infty$ while $n$ is fixed, have been studied by Hall, Marron and Neeman [8] and Ahn et al. [1]. They explored conditions which give the *geometric representation* of HDLSS data (i.e., modulo rotation, data tend to lie at vertices of a regular simplex) as well as strong inconsistency of eigenvectors. Strong inconsistency is also found in the context of $\frac{n}{d} \to c$, in the study of *phase transition*; see for example, Paul [16], Johnstone and Lu [12] and Baik, Ben Arous and Péché [2].

A reviewer pointed out a useful framework for organizing these variation is:

1. Classical: $d(n)/n \to 0$, as $n \to \infty$.
2. Random matrices: $d(n)/n \to c$, as $n \to \infty$.
3. HDLSS: $n$ fixed, with $d \to \infty$.

We view all of these as informative. Which is most informative will depend on the particular data analytic setting, in the same way that either the Normal or Poisson approximation can be "most informative" about the Binomial distribution.

In this paper, we focus only on the HDLSS case, and a broad and general set of conditions for consistency and strong inconsistency are provided. Section 2 develops conditions that guarantee the nonzero eigenvalues of the sample covariance matrix tend to an increasing constant, which are much more general than those of Hall, Marron and Neeman [8] and Ahn et al. [1]. This asymptotic behavior of the sample covariance matrix is the basis of the geometric representation of HDLSS data. Our result gives a broad new insight into this representation as discussed in Section 3. The central issue of consistency and strong inconsistency is developed in Section 4, as a series of theorems. For a fixed number $\kappa$, we assume the first $\kappa$ eigenvalues are much larger than the others. We show that when $\kappa = 1$, the first sample eigenvector is consistent and the others are strongly inconsistent. We also generalize to the $\kappa > 1$ case, featuring two different types of results (consistency and subspace consistency) according to the asymptotic behaviors of the first $\kappa$ eigenvalues. All results are combined and generalized in the main theorem (Theorem 2). Proofs of theorems are given in Section 5.

1.1. *General setting.* Suppose we have a $d \times n$ data matrix $X_{(d)} = [X_{1,(d)}, \ldots, X_{n,(d)}]$ with $d > n$, where the $d$-dimensional random vectors $X_{1,(d)}, \ldots, X_{n,(d)}$ are independent and identically distributed. We assume that each $X_{i,(d)}$ follows a multivariate distribution (which does not have to be Gaussian) with mean zero and covariance matrix $\Sigma_d$. Define the sphered data matrix $Z_{(d)} = \Lambda_d^{-1/2} U_d' X_{(d)}$. Then the components of the $d \times n$ matrix $Z_{(d)}$ have unit variances, and are uncorrelated with each other. We shall regulate the dependency (recall for non-Gaussian data, uncorrelated variables can still be dependent) of the random variables in $Z_{(d)}$ by a $\rho$-mixing condition. This allows serious weakening of the assumptions of Gaussianity while still enabling the law of large numbers that lie behind the geometric representation results of Hall, Marron and Neeman [8].

The concept of $\rho$-mixing was first developed by Kolmogorov and Rozanov [14]. See Bradley [5] for a clear and insightful discussion. For $-\infty \le J \le L \le \infty$, let $\mathcal{F}_J^L$ denote the $\sigma$-field of events generated by the random variables $(Z_i, J \le i \le L)$. For any $\sigma$-field $\mathcal{A}$, let $L_2(\mathcal{A})$ denote the space of square-integrable, $\mathcal{A}$ measurable (real-valued) random variables. For each $m \ge 1$, define the maximal correlation coefficient

$$\rho(m) := \sup |\operatorname{corr}(f,g)|, \qquad f \in L_2(\mathcal{F}_{-\infty}^j), g \in L_2(\mathcal{F}_{j+m}^\infty),$$

where sup is over all $f$, $g$ and $j \in \mathbf{Z}$. The sequence $\{Z_i\}$ is said to be $\rho$-mixing if $\rho(m) \to 0$ as $m \to \infty$.

While the concept of $\rho$-mixing is useful as a mild condition for the development of laws of large numbers, its formulation is critically dependent on the ordering of variables. For many interesting data types, such as microarray data, there is clear dependence but no natural ordering of the variables. Hence, we assume that there is some permutation of the data which is $\rho$-mixing. In particular, let $\{Z_{ij,(d)}\}_{i=1}^{d}$ be the components of the $j$th column vector of $Z_{(d)}$. We assume that for each $d$, there exists a permutation $\pi_d : \{1, \ldots, d\} \longmapsto \{1, \ldots, d\}$ so that the sequence $\{Z_{\pi_d(i)j,(d)} : i = 1, \ldots, d\}$ is $\rho$-mixing. This assumption makes the results invariant under a permutation of the variables.

In the following, all the quantities depend on $d$, but the subscript $d$ will be omitted for the sake of simplicity when it does not cause any confusion. The sample covariance matrix is defined as $S = n^{-1}XX'$. We do not subtract the sample mean vector because the population mean is assumed to be 0. Since the dimension of the sample covariance matrix $S$ grows, it is challenging to deal with $S$ directly. A useful approach is to work with the *dual* of $S$. The dual approach switches the role of columns and rows of the data matrix, by replacing $X$ by $X'$. The $n \times n$ *dual sample covariance matrix* is defined as $S_D = n^{-1}X'X$. An advantage of this dual approach is that $S_D$ and $S$ share nonzero eigenvalues. If we write $X$ as $U\Lambda^{1/2}Z$ and use the fact that $U$ is a unitary matrix,

$$(1.1) \qquad nS_D = (Z'\Lambda^{1/2}U')(U\Lambda^{1/2}Z) = Z'\Lambda Z = \sum_{i=1}^{d} \lambda_{i,d} z_i' z_i,$$

where the $z_i$'s, $i = 1, \ldots, d$, are the row vectors of the matrix $Z$. Note that $nS_D$ is commonly referred to as the *Gram matrix*, consisting of inner products between observations.

**2. HDLSS asymptotic behavior of the sample covariance matrix.** In this section, we investigate the behavior of the sample covariance matrix $S$ when $d \to \infty$ and $n$ is fixed. Under mild and broad conditions, the eigenvalues of $S$, or the dual $S_D$, behave asymptotically as if they are from the identity matrix. That is, the set of sample eigenvectors tends to be an arbitrary choice. This lies at the heart of the geometric representation results of Hall, Marron and Neeman [8] and Ahn et al. [1] which are studied more deeply in Section 3. We will see that this condition readily implies the strong inconsistency of sample eigenvectors; see Theorem 2.

The conditions for the theorem are conveniently formulated in terms of a measure of sphericity

$$\varepsilon \equiv \frac{\mathrm{tr}^2(\Sigma)}{d\,\mathrm{tr}(\Sigma^2)} = \frac{(\sum_{i=1}^{d} \lambda_{i,d})^2}{d\sum_{i=1}^{d} \lambda_{i,d}^2},$$

proposed and used by John [9, 10] as the basis of a hypothesis test for equality of eigenvalues. Note that these inequalities always hold:

$$\frac{1}{d} \leq \varepsilon \leq 1.$$

Also note that perfect sphericity of the distribution (i.e., equality of eigenvalues) occurs only when $\varepsilon = 1$. The other end of the $\varepsilon$ range is the most singular case where in the limit as the first eigenvalue dominates all others.

Ahn et al. [1] claimed that if $\varepsilon \gg \frac{1}{d}$, in the sense that $\varepsilon^{-1} = o(d)$, then the eigenvalues of $S_D$ tend to be identical in probability as $d \to \infty$. However, they needed an additional assumption (e.g., a Gaussian assumption on $X_{(d)}$) to have independence among components of $Z_{(d)}$, as described in Example 3.1. In this paper, we extend this result to the case of arbitrary distributions with dependency regulated by the $\rho$-mixing condition as in Section 1.1, which is much more general than either a Gaussian or an independence assumption. We also explore convergence in the almost sure sense with stronger assumptions. Our results use a measure of sphericity for part of the eigenvalues for conditions of a.s. convergence and also for later use in Section 4. In particular, define the measure of sphericity for $\{\lambda_{k,d}, \ldots, \lambda_{d,d}\}$ as

$$\varepsilon_k \equiv \frac{(\sum_{i=k}^d \lambda_{i,d})^2}{d \sum_{i=k}^d \lambda_{i,d}^2}.$$

For convenience, we name several assumptions used in this paper made about the measure of sphericity $\varepsilon$:

- *The $\varepsilon$-condition*: $\varepsilon \gg \frac{1}{d}$, that is,

$$(2.1) \qquad (d\varepsilon)^{-1} = \frac{\sum_{i=1}^d \lambda_{i,d}^2}{(\sum_{i=1}^d \lambda_{i,d})^2} \to 0 \qquad \text{as } d \to \infty.$$

- *The $\varepsilon_k$-condition*: $\varepsilon_k \gg \frac{1}{d}$, that is,

$$(2.2) \qquad (d\varepsilon_k)^{-1} = \frac{\sum_{i=k}^d \lambda_{i,d}^2}{(\sum_{i=k}^d \lambda_{i,d})^2} \to 0 \qquad \text{as } d \to \infty.$$

- *The strong $\varepsilon_k$-condition*: For some fixed $l \geq k$, $\varepsilon_l \gg \frac{1}{\sqrt{d}}$, that is,

$$(2.3) \qquad d^{-1/2} \varepsilon_l^{-1} = \frac{d^{1/2} \sum_{i=l}^d \lambda_{i,d}^2}{(\sum_{i=l}^d \lambda_{i,d})^2} \to 0 \qquad \text{as } d \to \infty.$$

REMARK. Note that the $\varepsilon_k$-condition is identical to the $\varepsilon$-condition when $k = 1$. Similarly, the strong $\varepsilon_k$-condition is also called *the strong $\varepsilon$-condition* when $k = 1$. The strong $\varepsilon_k$-condition is stronger than the $\varepsilon_k$ condition if the minimum of $l$'s which satisfy (2.3), $l_o$, is as small as $k$. But, if $l_o > k$, then this is not necessarily true. We will use the strong $\varepsilon_k$-condition combined with the $\varepsilon_k$-condition.

Note that the $\varepsilon$-condition is quite broad in the spectrum of possible values of $\varepsilon$: It only avoids the most singular case. The strong $\varepsilon$-condition further restricts $\varepsilon_l$ to essentially in the range $(\frac{1}{\sqrt{d}}, 1]$.

The following theorem states that if the (strong) $\varepsilon$-condition holds for $\Sigma_d$, then the sample eigenvalues behave as if they are from a scaled identity matrix. It uses the notation $I_n$ for the $n \times n$ identity matrix.

THEOREM 1. *For a fixed $n$, let $\Sigma_d = U_d \Lambda_d U_d'$, $d = n+1, n+2, \ldots$, be a sequence of covariance matrices. Let $X_{(d)}$ be a $d \times n$ data matrix from a $d$-variate distribution with mean zero and covariance matrix $\Sigma_d$. Let $S_d = \hat{U}_d \hat{\Lambda}_d \hat{U}_d'$ be the sample covariance matrix estimated from $X_{(d)}$ for each $d$ and let $S_{D,d}$ be its dual.*

*(1) Assume that the components of $Z_{(d)} = \Lambda_d^{-1/2} U_d' X_{(d)}$ have uniformly bounded fourth moments and are $\rho$-mixing under some permutation. If (2.1) holds, then*

$$(2.4) \qquad c_d^{-1} S_{D,d} \longrightarrow I_n,$$

*in probability as $d \to \infty$, where $c_d = n^{-1} \sum_{i=1}^{d} \lambda_{i,d}$.*

*(2) Assume that the components of $Z_{(d)} = \Lambda_d^{-1/2} U_d' X_{(d)}$ have uniformly bounded eighth moments and are independent to each other. If both (2.1) and (2.3) hold, then $c_d^{-1} S_{D,d} \to I_n$ almost surely as $d \to \infty$.*

The (strong) $\varepsilon$-condition holds for quite general settings. The strong $\varepsilon$-condition combined with the $\varepsilon$-condition holds under:

(a) Null case: All eigenvalues are the same.
(b) Mild spiked model: The first $m$ eigenvalues are moderately larger than the others, for example, $\lambda_{1,d} = \cdots = \lambda_{m,d} = C_1 \cdot d^\alpha$ and $\lambda_{m+1,d} = \cdots = \lambda_{d,d} = C_2$, where $m < d$, $\alpha < 1$ and $C_1, C_2 > 0$.

The $\varepsilon$-condition fails when:

(c) Singular case: Only the first few eigenvalues are nonzero.
(d) Exponential decrease: $\lambda_{i,d} = c^{-i}$ for some $c > 1$.
(e) Sharp spiked model: The first $m$ eigenvalues are much larger than the others. One example is the same as (b), but $\alpha \geq 1$.

The polynomially decreasing case, $\lambda_{i,d} = i^{-\beta}$, is interesting because it depends on the power $\beta$:

(f-1) The strong $\varepsilon$-condition holds when $0 \leq \beta < \frac{3}{4}$.
(f-2) The $\varepsilon$-condition holds, but the strong $\varepsilon$-condition fails when $\frac{3}{4} \leq \beta \leq 1$.
(f-3) The $\varepsilon$-condition fails when $\beta > 1$.

Another family of examples that includes all three cases is the spiked model with the number of spikes increasing, for example, $\lambda_{1,d} = \cdots = \lambda_{m,d} = C_1 \cdot d^\alpha$ and $\lambda_{m+1,d} = \cdots = \lambda_{d,d} = C_2$, where $m = \lfloor d^\beta \rfloor$, $0 < \beta < 1$ and $C_1, C_2 > 0$:

(g-1) The strong $\varepsilon$-condition holds when $0 \leq 2\alpha + \beta < \frac{3}{2}$.
(g-2) The $\varepsilon$-condition holds but the strong $\varepsilon$-condition fails when $\frac{3}{2} \leq 2\alpha + \beta < 2$.
(g-3) The $\varepsilon$-condition fails when $2\alpha + \beta \geq 2$.

**3. Geometric representation of HDLSS data.** Suppose $X \sim \mathcal{N}_d(0, I_d)$. When the dimension $d$ is small, most of the mass of the data lies near origin. However, with a large $d$, Hall, Marron and Neeman [8] showed that Euclidean distance of $X$ to the origin is described as

$$\|X\| = \sqrt{d} + o_p(\sqrt{d}). \tag{3.1}$$

Moreover, the distance between two samples is also rather deterministic, that is,

$$\|X_1 - X_2\| = \sqrt{2d} + o_p(\sqrt{d}). \tag{3.2}$$

These results can be derived by the law of large numbers. Hall, Marron and Neeman [8] generalized those results under the assumptions that $d^{-1} \sum_{i=1}^d \text{Var}(X_i) \to 1$ and $\{X_i\}$ is $\rho$-mixing.

Application of part (1) of Theorem 1 generalizes these results. Let $X_{1,(d)}$, $X_{2,(d)}$ be two samples that satisfy the assumptions of Theorem 1 part (1). Assume without loss of generality that $\lim_{d \to \infty} d^{-1} \sum_{i=1}^d \lambda_{i,d} = 1$. The scaled squared distance between two data points is

$$\frac{\|X_{1,(d)} - X_{2,(d)}\|^2}{\sum_{i=1}^d \lambda_{i,d}} = \sum_{i=1}^d \tilde{\lambda}_{i,d} z_{i1}^2 + \sum_{i=1}^d \tilde{\lambda}_{i,d} z_{i2}^2 - 2 \sum_{i=1}^d \tilde{\lambda}_{i,d} z_{i1} z_{i2},$$

where $\tilde{\lambda}_{i,d} = \frac{\lambda_{i,d}}{\sum_{i=1}^d \lambda_{i,d}}$. Note that by (1.1), the first two terms are diagonal elements of $c_d^{-1} S_{D,d}$ in Theorem 1 and the third term is an off-diagonal element. Since $c_d^{-1} S_{D,d} \to I_n$, we have (3.2). (3.1) is derived similarly.

REMARK. If $\lim_{d \to \infty} d^{-1} \sum_{i=1}^d \lambda_{i,d} = 1$, then the conclusion (2.4) of Theorem 1 part (1) holds if and only if the representations (3.1) and (3.2) hold under the same assumptions in the theorem.

In this representation, the $\rho$-mixing assumption plays a very important role. The following example, due to John Kent, shows that some type of mixing condition is important.

EXAMPLE 3.1 (Strong dependency via a scale mixture of Gaussian). Let $X = Y_1 U + \sigma Y_2 (1 - U)$, where $Y_1, Y_2$ are two independent $\mathcal{N}_d(0, I_d)$ random variables, $U = 0$ or $1$ with probability $\frac{1}{2}$ and independent of $Y_1, Y_2$, and $\sigma > 1$. Then,

$$\|X\| = \begin{cases} d^{1/2} + O_p(1), & \text{w.p. } \frac{1}{2}, \\ \sigma d^{1/2} + O_p(1), & \text{w.p. } \frac{1}{2}. \end{cases}$$

Thus, (3.1) does not hold. Note that since $\mathrm{Cov}(X) = \frac{1+\sigma^2}{2} I_d$, the $\varepsilon$-condition holds and the variables are uncorrelated. However, there is strong dependency, i.e., $\mathrm{Cov}(z_i^2, z_j^2) = (\frac{1+\sigma^2}{2})^{-2} \mathrm{Cov}(x_i^2, x_j^2) = (\frac{1-\sigma^2}{1+\sigma^2})^2$ for all $i \neq j$ which implies that $\rho(m) > c$ for some $c > 0$, for all $m$. Thus, the $\rho$-mixing condition does not hold for all permutation. Note that, however, under Gaussian assumption, given any covariance matrix $\Sigma$, $Z = \Sigma^{-1/2} X$ has independent components.

Note that in the case $X = (X_1, \ldots, X_d)$ is a sequence of i.i.d. random variables, the results (3.1) and (3.2) can be considerably strengthened to $\|X\| = \sqrt{d} + O_p(1)$, and $\|X_1 - X_2\| = \sqrt{2d} + O_p(1)$. The following example shows that strong results are beyond the reach of reasonable assumption.

EXAMPLE 3.2 (Varying sphericity). Let $X \sim \mathcal{N}_d(0, \Sigma_d)$, where $\Sigma_d = \mathrm{diag}(d^\alpha, 1, \ldots, 1)$ and $\alpha \in (0, 1)$. Define $Z = \Sigma_d^{-1/2} X$. Then the components of $Z$, $z_i$'s, are independent standard Gaussian random variables. We get $\|X\|^2 = d^\alpha z_1^2 + \sum_{i=2}^{d} z_i^2$. Now for $0 < \alpha < \frac{1}{2}$, $d^{-1/2}(\|X\|^2 - d) \Rightarrow \mathcal{N}(0, 1)$ and for $\frac{1}{2} < \alpha < 1$, $d^{-\alpha}(\|X\|^2 - d) \Rightarrow z_1^2$, where $\Rightarrow$ denotes convergence in distribution. Thus, by the delta-method, we get

$$\|X\| = \begin{cases} \sqrt{d} + O_p(1), & \text{if } 0 < \alpha < \frac{1}{2}, \\ \sqrt{d} + O_p(d^{\alpha - 1/2}), & \text{if } \frac{1}{2} < \alpha < 1. \end{cases}$$

In both cases, the representation (3.1) holds.

**4. Consistency and strong inconsistency of PC directions.** In this section, conditions for consistency or strong inconsistency of the sample PC direction vectors are investigated in the general setting of Section 1.1. The generic eigen-structure of the covariance matrix that we assume is the following. For a fixed number $\kappa$, we assume the first $\kappa$ eigenvalues are much larger than others. (The precise meaning of *large* will be addressed shortly.)

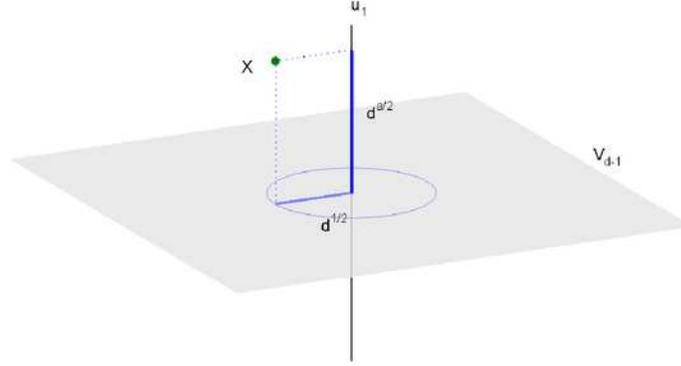

Fig. 2. *Projection of a d-dimensional random variable $X$ onto $u_1$ and $V_{d-1}$. If $\alpha > 1$, then the subspace $V_{d-1}$ becomes negligible compared to $u_1$ when $d \to \infty$.*

The rest of eigenvalues are assumed to satisfy the $\varepsilon$-condition, which is very broad in the range of sphericity. We begin with the case $\kappa = 1$ and generalize the result for $\kappa > 1$ in two distinct ways. The main theorem (Theorem 2) contains and combines those previous results and also embraces various cases according to the magnitude of the first $\kappa$ eigenvalues. We also investigate the sufficient conditions for a stronger result, that is, almost sure convergence, which involves use of the strong $\varepsilon$-condition.

4.1. *Criteria for consistency or strong inconsistency of the first PC direction.* Consider the simplest case that only the first PC direction of $S$ is of interest. Section 3 gives some preliminary indication of this. As an illustration, consider a spiked model as in Example 3.2 but now let $\alpha > 1$. Let $\{u_i\}$ be the set of eigenvectors of $\Sigma_d$ and $V_{d-1}$ be the subspace of all eigenvectors except the first one. Then the projection of $X$ onto $u_1$ has a norm $\|\operatorname{Proj}_{u_1} X\| = \|X_1\| = O_p(d^{\alpha/2})$. The projection of $X$ onto $V_{d-1}$ has a norm $\sqrt{d} + o_p(\sqrt{d})$ by (3.1). Thus, when $\alpha > 1$, if we scale the whole data space $\mathbf{R}^d$ by dividing by $d^{\alpha/2}$, then $\operatorname{Proj}_{V_{d-1}} X$ becomes negligible compared to $\operatorname{Proj}_{u_1} X$ (see Figure 2). Thus, for a large $d$, $\Sigma_d \approx \lambda_1 u_1 u_1'$ and the variation of $X$ is mostly along $u_1$. Therefore, the sample eigenvector corresponding to the largest eigenvalue, $\hat{u}_1$, will be similar to $u_1$.

To generalize this, suppose the $\varepsilon_2$ condition holds. The following proposition states that under the general setting in Section 1.1, the first sample eigenvector $\hat{u}_1$ converges to its population counterpart $u_1$ (consistency) or tends to be perpendicular to $u_1$ (strong inconsistency) according to the magnitude of the first eigenvalue $\lambda_1$, while all the other sample eigenvectors are strongly inconsistent regardless of the magnitude $\lambda_1$.

PROPOSITION 1. *For a fixed $n$, let $\Sigma_d = U_d \Lambda_d U_d'$, $d = n+1, n+2, \ldots$, be a sequence of covariance matrices. Let $X_{(d)}$ be a $d \times n$ data matrix from* a $d$-variate distribution with mean zero and covariance matrix $\Sigma_d$. Let $S_d = \hat{U}_d \hat{\Lambda}_d \hat{U}'_d$ be the sample covariance matrix estimated from $X_{(d)}$ for each $d$. Assume the following:

(a) The components of $Z_{(d)} = \Lambda_d^{-1/2} U'_d X_{(d)}$ have uniformly bounded fourth moments and are $\rho$-mixing for some permutation.

For an $\alpha_1 > 0$,

(b) $\frac{\lambda_{1,d}}{d^{\alpha_1}} \longrightarrow c_1$ for some $c_1 > 0$.
(c) The $\varepsilon_2$-condition holds and $\sum_{i=2}^{d} \lambda_{i,d} = O(d)$.

If $\alpha_1 > 1$, then the first sample eigenvector is consistent and the others are strongly inconsistent in the sense that

$$\text{Angle}(\hat{u}_1, u_1) \xrightarrow{p} 0 \quad \text{as } d \to \infty,$$

$$\text{Angle}(\hat{u}_i, u_i) \xrightarrow{p} \frac{\pi}{2} \quad \text{as } d \to \infty \ \forall i = 2, \ldots, n.$$

If $\alpha_1 \in (0, 1)$, then all sample eigenvectors are strongly inconsistent, i.e.,

$$\text{Angle}(\hat{u}_i, u_i) \xrightarrow{p} \frac{\pi}{2} \quad \text{as } d \to \infty \ \forall i = 1, \ldots, n.$$

Note that the gap between consistency and strong inconsistency is very thin, i.e., if we avoid $\alpha_1 = 1$, then we have either consistency or strong inconsistency. Thus in the HDLSS context, asymptotic behavior of PC directions is mostly captured by consistency and strong inconsistency. Now it makes sense to say $\lambda_1$ is much larger than the others when $\alpha_1 > 1$, which results in consistency. Also note that if $\alpha_1 < 1$, then the $\varepsilon$-condition holds, which is in fact the condition for Theorem 1.

4.2. *Generalizations.* In this section, we generalize Proposition 1 to the case that multiple eigenvalues are much larger than the others. This leads to two different types of result.

First is the case that the first $p$ eigenvectors are each consistent. Consider a covariance structure with multiple spikes, that is, $p$ eigenvalues, $p > 1$, which are much larger than the others. In order to have consistency of the first $p$ eigenvectors, we require that each of $p$ eigenvalues has a distinct order of magnitude, for example, $\lambda_{1,d} = d^3$, $\lambda_{2,d} = d^2$ and sum of the rest is order of $d$.

PROPOSITION 2. *For a fixed $n$, let $\Sigma_d$, $X_{(d)}$, and $S_d$ be as before. Assume (a) of Proposition 1. Let $\alpha_1 > \alpha_2 > \cdots > \alpha_p > 1$ for some $p < n$. Suppose the following conditions hold:*

(b) $\frac{\lambda_{i,d}}{d^{\alpha_i}} \longrightarrow c_i$ for some $c_i > 0 \ \forall i = 1, \ldots, p$.

(c) *The $\varepsilon_{p+1}$-condition holds and $\sum_{i=p+1}^{d} \lambda_{i,d} = O(d)$.*

*Then, the first $p$ sample eigenvectors are consistent and the others are strongly inconsistent in the sense that*

$$\mathrm{Angle}(\hat{u}_i, u_i) \xrightarrow{p} 0 \qquad as\ d \to \infty\ \forall i = 1, \ldots, p,$$

$$\mathrm{Angle}(\hat{u}_i, u_i) \xrightarrow{p} \frac{\pi}{2} \qquad as\ d \to \infty\ \forall i = p+1, \ldots, n.$$

Consider now a distribution having a covariance structure with multiple spikes as before. Let $k$ be the number of spikes. An interesting phenomenon happens when the first $k$ eigenvalues are of the same order of magnitude, that is, $\lim_{d \to \infty} \frac{\lambda_{1,d}}{\lambda_{k,d}} = c > 1$ for some constant $c$. Then the first $k$ sample eigenvectors are neither consistent nor strongly inconsistent. However, all of those random directions converge to the subspace spanned by the first $k$ population eigenvectors. Essentially, when eigenvalues are of the same order, the eigen-directions can not be separated but are subspace consistent with the proper subspace.

PROPOSITION 3. *For a fixed $n$, let $\Sigma_d$, $X_{(d)}$, and $S_d$ be as before. Assume* (a) *of Proposition 1. Let $\alpha_1 > 1$ and $k < n$. Suppose the following conditions hold:*

(b) $\frac{\lambda_{i,d}}{d^{\alpha_1}} \longrightarrow c_i$ *for some $c_i > 0$ $\forall i = 1, \ldots, k$.*
(c) *The $\varepsilon_{k+1}$-condition holds and $\sum_{i=k+1}^{d} \lambda_{i,d} = O(d)$.*

*Then the first $k$ sample eigenvectors are subspace-consistent with the subspace spanned by the first $k$ population eigenvectors, and the others are strongly inconsistent in the sense that*

$$\mathrm{Angle}(\hat{u}_i, \mathrm{span}\{u_1, \ldots, u_k\}) \xrightarrow{p} 0 \qquad as\ d \to \infty\ \forall i = 1, \ldots, k,$$

$$\mathrm{Angle}(\hat{u}_i, u_i) \xrightarrow{p} \frac{\pi}{2} \qquad as\ d \to \infty\ \forall i = k+1, \ldots, n.$$

4.3. *Main theorem.* Propositions 1–3 are combined and generalized in the main theorem. Consider $p$ groups of eigenvalues, which grow at the same rate within each group as in Proposition 3. Each group has a finite number of eigenvalues and the number of eigenvalues in all groups, $\kappa$, does not exceed $n$. Also similar to Proposition 2, let the orders of magnitude of the $p$ groups be different to each other. We require that the $\varepsilon_{\kappa+1}$-condition holds. The following theorem states that a sample eigenvector of a group converges to the subspace of population eigenvectors of the group.

THEOREM 2 (Main theorem). *For a fixed $n$, let $\Sigma_d$, $X_{(d)}$, and $S_d$ be as before. Assume* (a) *of Proposition 1. Let $\alpha_1, \ldots, \alpha_p$ be such that $\alpha_1 > \alpha_2 >$* $\cdots > \alpha_p > 1$ *for some* $p < n$. *Let* $k_1, \ldots, k_p$ *be nonnegative integers such that* $\sum_{j=1}^p k_j \doteq \kappa < n$. *Let* $k_0 = 0$ *and* $k_{p+1} = d - \kappa$. *Let* $J_1, \ldots, J_{p+1}$ *be sets of indices such that*

$$J_l = \left\{1 + \sum_{j=0}^{l-1} k_j, 2 + \sum_{j=0}^{l-1} k_j, \ldots, k_l + \sum_{j=0}^{l-1} k_j\right\}, \qquad l = 1, \ldots, p+1.$$

*Suppose the following conditions hold:*

(b) $\frac{\lambda_{i,d}}{d^{\alpha_l}} \longrightarrow c_i$ *for some* $c_i > 0, \forall i \in J_l, \forall l = 1, \ldots, p$.
(c) *The* $\varepsilon_{\kappa+1}$-*condition holds and* $\sum_{i \in J_{p+1}} \lambda_{i,d} = O(d)$.

*Then the sample eigenvectors whose label is in the group* $J_l$, *for* $l = 1, \ldots, p$, *are subspace-consistent with the space spanned by the population eigenvectors whose labels are in* $J_l$ *and the others are strongly inconsistent in the sense that*

$$(4.1) \quad \mathrm{Angle}(\hat{u}_i, \mathrm{span}\{u_j : j \in J_l\}) \xrightarrow{p} 0 \qquad as\ d \to \infty\ \forall i \in J_l, \forall l = 1, \ldots, p,$$

*and*

$$(4.2) \qquad \mathrm{Angle}(\hat{u}_i, u_i) \xrightarrow{p} \frac{\pi}{2} \qquad as\ d \to \infty\ \forall i = \kappa+1, \ldots, n.$$

REMARK. If the cardinality of $J_l$, $k_l$, is 1, then (4.1) implies $\hat{u}_i$ is consistent for $i \in J_l$.

REMARK. The strongly inconsistent eigenvectors whose labels are in $J_{p+1}$ can be considered to be subspace-consistent. Let $\Gamma_d$ be the subspace spanned by the population eigenvectors whose labels are in $J_{p+1}$ for each $d$, i.e. $\Gamma_d = \mathrm{span}\{u_j : j \in J_{p+1}\} = \mathrm{span}\{u_{\kappa+1}, \ldots, u_d\}$. Then

$$\mathrm{Angle}(\hat{u}_{i,d}, \Gamma_d) \xrightarrow{p} 0 \qquad as\ d \to \infty$$

for all $i \in J_{p+1}$.

Note that the formulation of the theorem is similar to the spiked covariance model but much more general. The uniform assumption on the underlying eigenvalues, that is, $\lambda_i = 1$ for all $i > \kappa$, is relaxed to the $\varepsilon$-condition. We also have catalogued a large collection of specific results according to the various sizes of spikes.

These results are now illustrated for some classes of covariance matrices that are of special interest. These covariance matrices are easily represented in *factor form*, that is, in terms of $F_d = \Sigma_d^{1/2}$.

EXAMPLE 4.1. Consider a series of covariance matrices $\{\Sigma_d\}_d$. Let $\Sigma_d = F_d F_d'$, where $F_d$ is a $d \times d$ symmetric matrix such that

$$F_d = (1 - \rho_d)I_d + \rho_d J_d = \begin{pmatrix} 1 & \rho_d & \cdots & \rho_d \\ \rho_d & 1 & \ddots & \vdots \\ \vdots & \ddots & \ddots & \rho_d \\ \rho_d & \cdots & \rho_d & 1 \end{pmatrix},$$

where $J_d$ is the $d \times d$ matrix of ones and $\rho_d \in (0,1)$ depends on $d$. The eigenvalues of $\Sigma_d$ are $\lambda_{1,d} = (d\rho_d + 1 - \rho_d)^2, \lambda_{2,d} = \cdots = \lambda_{d,d} = (1 - \rho_d)^2$. Note that this is a simple and natural probabilistic mechanism that generates eigenvalues where the first is order of magnitude larger than the rest (our fundamental assumption). The first eigenvector is $u_1 = \frac{1}{\sqrt{d}}(1,1,\ldots,1)'$, while $\{u_2,\ldots,u_d\}$ are any orthogonal sets of direction vectors perpendicular to $u_1$. Note that $\sum_{i=2}^d \lambda_{i,d} = d(1-\rho_d)^2 = O(d)$ and the $\varepsilon_2$-condition holds. Let $X_d \sim \mathcal{N}_d(0, \Sigma_d)$. By Theorem 2, if $\rho_d \in (0,1)$ is a fixed constant or decreases to 0 slowly so that $\rho_d \gg d^{-1/2}$, then the first PC direction $\hat{u}_1$ is consistent. Else if $\rho_d$ decreases to 0 so quickly that $\rho_d \ll d^{-1/2}$, then $\hat{u}_1$ is strongly inconsistent. In both cases, all the other sample PC directions are strongly inconsistent.

EXAMPLE 4.2. Consider now a $2d \times 2d$ covariance matrix $\Sigma_d = F_d F_d'$, where $F_d$ is a block diagonal matrix, such that

$$F_d = \begin{pmatrix} F_{1,d} & O \\ O & F_{2,d} \end{pmatrix},$$

where $F_{1,d} = (1-\rho_{1,d})I_d + \rho_{1,d}J_d$ and $F_{2,d} = (1-\rho_{2,d})I_d + \rho_{2,d}J_d$. Suppose $0 < \rho_{2,d} \leq \rho_{1,d} < 1$. Note that $\lambda_{1,d} = (d\rho_{1,d} + 1 - \rho_{1,d})^2$, $\lambda_{2,d} = (d\rho_{2,d} + 1 - \rho_{2,d})^2$ and the $\varepsilon_3$-condition holds. Let $X_{2d} \sim \mathcal{N}_{2d}(0, \Sigma_d)$. Application of Theorem 2 for various conditions on $\rho_{1,d}, \rho_{2,d}$ is summarized as follows. Denote, for two nonincreasing sequences $\mu_d, \nu_d \in (0,1)$, $\mu_d \gg \nu_d$ for $\nu_d = o(\mu_d)$ and $\mu_d \succeq \nu_d$ for $\lim_{d \to \infty} \frac{\mu_d}{\nu_d} = c \in [1, \infty)$:

1. $\rho_{1,d} \gg \rho_{2,d} \gg d^{-1/2}$: Both $\hat{u}_1, \hat{u}_2$ consistent.
2. $\rho_{1,d} \succeq \rho_{2,d} \gg d^{-1/2}$: Both $\hat{u}_1, \hat{u}_2$ subspace-consistent to span$\{u_1, u_2\}$.
3. $\rho_{1,d} \gg d^{-1/2} \gg \rho_{2,d}$: $\hat{u}_1$ consistent, $\hat{u}_2$ strongly inconsistent.
4. $d^{-1/2} \gg \rho_{1,d} \gg \rho_{2,d}$: Both $\hat{u}_1, \hat{u}_2$ strongly inconsistent.

4.4. *Corollaries to the main theorem.* The result can be extended for special cases.

First of all, consider constructing $X_{(d)}$ from $Z_d$ by $X_{(d)} \equiv U_d \Lambda_d^{1/2} Z_d$ where $Z_d$ is a truncated set from an infinite sequence of independent random variables with mean zero and variance 1. This assumption makes it possible to have convergence in the almost sure sense. This is mainly because the triangular array $\{Z_{1i,(d)}\}_{i,d}$ becomes the single sequence $\{Z_{1i}\}_i$.

COROLLARY 1. *Suppose all the assumptions in Theorem 2, with the assumption* (a) *replaced by the following:*

(a′) *The components of $Z_{(d)} = \Lambda_d^{-1/2} U_d' X_{(d)}$ have uniformly bounded eighth moments and are independent to each other. Let $Z_{1i,(d)} \equiv Z_{1i}$ for all $i, d$.*

*If the strong $\varepsilon_{\kappa+1}$-condition (2.3) holds, then the mode of convergence of (4.1) and (4.2) is almost sure.*

Second, consider the case that both $d, n$ tend to infinity. Under the setting of Theorem 2, we can separate PC directions better when the eigenvalues are distinct. When $d \to \infty$, we have subspace consistency of $\hat{u}_i$ with the proper subspace, which includes $u_i$. Now letting $n \to \infty$ makes it possible for $\hat{u}_i$ to be consistent.

COROLLARY 2. *Let $\Sigma_d$, $X_{(d)}$ and $S_d$ be as before. Under the assumptions* (a), (b) *and* (c) *in Theorem 2, assume further for* (b) *that the first $\kappa$ eigenvalues are distinct, that is, $c_i > c_j$ for $i > j$ and $i, j \in J_l$ for $l = 1, \ldots, p$. Then for all $i \leq \kappa$,*

(4.3) $$\text{Angle}(\hat{u}_i, u_i) \xrightarrow{p} 0 \quad \text{as } d \to \infty,\, n \to \infty,$$

*where the limits are applied successively.*

*If the assumption* (a) *is replaced by the assumption* (a′) *of Corollary 1, then the mode of convergence of (4.3) is almost sure.*

This corollary can be viewed as the case when $d, n$ tend to infinity together, but $d$ increases at a much faster rate than $n$, that is, $d \gg n$. When $n$ also increases in the particular setting of the corollary, the sample eigenvectors, which were only subspace-consistent in the $d \to \infty$ case, tend to be distinguishable and each of the eigenvectors is consistent. We conjecture that the inconsistent sample eigenvalues are still strongly inconsistent when $d, n \to \infty$ and $d \gg n$.

4.5. *Limiting distributions of corresponding eigenvalues.* The study of asymptotic behavior of the sample eigenvalues is an important part in the proof of Theorem 2, and also could be of independent interest. The following lemma states that the large sample eigenvalues increase at the same speed as their population counterpart and the relatively small eigenvalues tend to be of order of $d$ as $d$ tends to infinity. Let $\varphi_i(A)$ denote the $i$th largest eigenvalue of the symmetric matrix $A$ and $\varphi_{i,l}(A) = \varphi_{i^*}(A)$ where $i^* = i - \sum_{j=1}^{l-1} k_j$.

16 S. JUNG AND J. S. MARRONLEMMA 1. *If the assumptions of Theorem 2 hold, and let $Z_l$ be a $k_l \times n$ matrix from blocks of $Z$ as defined in (5.2), then*

$$\hat{\lambda}_i/d^{\alpha_l} \Longrightarrow \eta_i \qquad \text{as } d \to \infty \text{ if } i \in J_l \ \forall l = 1, \ldots, p,$$

$$\hat{\lambda}_i/d \xrightarrow{p} K \qquad \text{as } d \to \infty \text{ if } i = \kappa + 1, \ldots, n,$$

*where each $\eta_i$ is a random variable whose support is $(0, \infty)$ almost surely and indeed $\eta_i = \varphi_{i,l}(n^{-1} C_l^{1/2} Z_l Z_l' C^{1/2})$ for each $i \in J_l$, where $C_l = \mathrm{diag}\{c_j : j \in J_l\}$ and $K = \lim_{d \to \infty} (dn)^{-1} \sum_{i \in J_{p+1}} \lambda_{i,d}$.*

If the data matrix $X_{(d)}$ is Gaussian, then the first $\kappa$ sample eigenvalues converge in distribution to some quantities, which have known distributions.

COROLLARY 3. *Under all the assumptions of Theorem 2, assume further that $X_{(d)} \sim \mathcal{N}_d(0, \Sigma_d)$ for each $d$. Then, for $i \in J_l$, $l = 1, \ldots, p$,*

$$\frac{\hat{\lambda}_i}{d^{\alpha_l}} \Longrightarrow \varphi_{i,l}(n^{-1} \mathcal{W}_{k_l}(n, C_l)) \qquad \text{as } d \to \infty,$$

*where $\mathcal{W}_{k_l}(n, C_l)$ denotes a $k_l \times k_l$ random matrix distributed as the Wishart distribution with degree of freedom $n$ and covariance $C_l$.*

*If $k_l = 1$ for some $l$, then for $i \in J_l$*

$$\frac{\hat{\lambda}_i}{\lambda_i} \Longrightarrow \frac{\chi_n^2}{n} \qquad \text{as } d \to \infty,$$

*where $\chi_n^2$ denotes a random variable distributed as the $\chi^2$ distribution with degree of freedom $n$.*

This generalizes the results in Section 4.2 of Ahn et al. [1].

**5. Proofs.**

5.1. *Proof of Theorem 1.* First, we give the proof of part (1). By (1.1), the $m$th diagonal entry of $nS_D$ can be expressed as $\sum_{i=1}^{d} \lambda_{i,d} z_{im,d}^2$ where $z_{im,d}$ is the $(i, m)$th entry of the matrix $Z_{(d)}$. Define the relative eigenvalues $\tilde{\lambda}_{i,d}$ as $\tilde{\lambda}_{i,d} \equiv \frac{\lambda_{i,d}}{\sum_{i=1}^{d} \lambda_{i,d}}$. Let $\pi_d$ denote the given permutation for each $d$ and let $Y_i = z_{\pi_d(i)m,d}^2 - 1$. Then the $Y_i$'s are $\rho$-mixing, $\mathrm{E}(Y_i) = 0$ and $\mathrm{E}(Y_i^2) \leq B$ for all $i$ for some $B < \infty$. Let $\rho(m) = \sup |\mathrm{corr}(Y_i, Y_{i+m})|$ where the sup is over all $i$. We shall use the following lemma.

LEMMA 2. *For any permutation $\pi_d^*$,*

$$\lim_{d \to \infty} \sum_{i=1}^{d} \tilde{\lambda}_{\pi_d^*(i),d} \rho(i) = 0.$$

PROOF. For any $\delta > 0$, since $\lim_{i \to \infty} \rho(i) = 0$, we can choose $N$ such that $\rho(i) < \delta/2$ for all $i > N$. Since $\lim_{d \to \infty} \sum_{i=1}^{d} \tilde{\lambda}^2_{\pi^*_d(i),d} = 0$, we get $\lim_{d \to \infty} \sum_{i=1}^{N} \tilde{\lambda}_{\pi^*_d(i),d} = 0$. Thus, we can choose $d_0$ satisfying $\sum_{i=1}^{N} \tilde{\lambda}_{\pi^*_d(i),d} < \frac{\delta}{2}$ for all $d > d_0$. With the fact $\sum_{i=1}^{d} \tilde{\lambda}_{i,d} = 1$ for all $d$ and $\rho(i) < 1$, we get for all $d > d_0$,

$$\sum_{i=1}^{d} \tilde{\lambda}_{\pi^*_d(i),d} \rho(i) = \sum_{i=1}^{N} \tilde{\lambda}_{\pi^*_d(i),d} \rho(i) + \sum_{i=N+1}^{d} \tilde{\lambda}_{\pi^*_d(i),d} \rho(i) < \delta. \qquad \Box$$

Now let $\pi_d^{-1}$ be the inverse permutation of $\pi_d$. Then by Lemma 2 and the $\varepsilon$-condition, there exists a permutation $\pi^*_d$ such that

$$\mathrm{E}\left(\sum_{i=1}^{d} \tilde{\lambda}_{\pi_d^{-1}(i),d} Y_i\right)^2 = \sum_{i=1}^{d} \tilde{\lambda}^2_{\pi_d^{-1}(i),d} \mathrm{E} Y_i^2 + 2\sum_{i=1}^{d} \tilde{\lambda}_{\pi_d^{-1}(i),d} \sum_{j=i+1}^{d} \tilde{\lambda}_{\pi_d^{-1}(j),d} \mathrm{E} Y_i Y_j$$

$$\leq \sum_{i=1}^{d} \tilde{\lambda}^2_{i,d} B + 2\sum_{i=1}^{d} \tilde{\lambda}_{i,d} \sum_{j=1}^{d} \tilde{\lambda}_{\pi^*_d(j),d} \rho(j) B^2 \to 0,$$

as $d \to \infty$. Then Chebyshev's inequality gives us, for any $\tau > 0$,

$$P\left[\left|\sum_{i=1}^{d} \tilde{\lambda}_{i,d} z_{im}^2 - 1\right| > \tau\right] \leq \frac{\mathrm{E}(\sum_{i=1}^{d} \tilde{\lambda}_{\pi_d^{-1}(i),d} Y_i)^2}{\tau^2} \to 0,$$

as $d \to \infty$. Thus, we conclude that the diagonal elements of $nS_D$ converge to 1 in probability.

The off-diagonal elements of $nS_D$ can be expressed as $\sum_{i=1}^{d} \lambda_{i,d} z_{im} z_{il}$. Similar arguments to those used in the diagonal case, together with the fact that $z_{im}$ and $z_{il}$ are independent, gives that

$$\mathrm{E}\left(\sum_{i=1}^{d} \tilde{\lambda}_{i,d} z_{im} z_{il}\right)^2 \leq \sum_{i=1}^{d} \tilde{\lambda}^2_{i,d} + 2\sum_{i=1}^{d} \tilde{\lambda}_{i,d} \sum_{j=i+1}^{d} \tilde{\lambda}_{\pi_d^{-1}(j),d} \rho^2(j-i) \to 0,$$

as $d \to \infty$. Thus, by Chebyshev's inequality, the off-diagonal elements of $nS_D$ converge to 0 in probability.

Now, we give the proof for part (2). We begin with the $m$th diagonal entry of $nS_D$, $\sum_{i=1}^{d} \lambda_{i,d} z_{im}^2$. Note that since $\sum_{i=1}^{k-1} \tilde{\lambda}_{i,d} \to 0$ by the $\varepsilon$-condition, we assume $k=1$ in (2.3) without loss of generality.

Let $Y_i = z_{im}^2 - 1$. Note that the $Y_i$'s are independent, $\mathrm{E}(Y_i) = 0$ and $\mathrm{E}(Y_i^4) \leq B$ for all $i$ for some $B < \infty$. Now

$$(5.1) \qquad \mathrm{E}\left(\sum_{i=1}^{d} \tilde{\lambda}_{i,d} Y_i\right)^4 = \mathrm{E} \sum_{i,j,k,l=1}^{d} \tilde{\lambda}_{i,d} \tilde{\lambda}_{j,d} \tilde{\lambda}_{k,d} \tilde{\lambda}_{l,d} Y_i Y_j Y_k Y_l.$$

Note that terms in the sum of the form $EY_iY_jY_kY_l$, $EY_i^2Y_jY_k$ and $EY_i^3Y_j$ are 0 if $i,j,k,l$ are distinct. The only terms that do not vanish are those of the form $EY_i^4$, $EY_i^2Y_j^2$, both of which are bounded by $B$. Note that $\tilde{\lambda}_{i,d}^2$'s are nonnegative, and hence the sum of squares is less than the square of sum, we have $\sum_{i=1}^d \tilde{\lambda}_{i,d}^4 \leq (\sum_{i=1}^d \tilde{\lambda}_{i,d}^2)^2$. Also note that by the strong $\varepsilon$-condition, $\sum_{i=1}^d \tilde{\lambda}_{i,d}^2 = (d\varepsilon)^{-1} = o(d^{-1/2})$. Thus, (5.1) is bounded as

$$E\left(\sum_{i=1}^d \tilde{\lambda}_{i,d}Y_i\right)^4 \leq \sum_{i=1}^d \tilde{\lambda}_{i,d}^4 B + \sum_{i=j\neq k=l} \tilde{\lambda}_{i,d}^2 \tilde{\lambda}_{k,d}^2 B$$

$$\leq \left(\sum_{i=1}^d \tilde{\lambda}_{i,d}^2\right)^2 B + \binom{4}{2}\left(\sum_{i=1}^d \tilde{\lambda}_{i,d}^2\right)^2 B$$

$$= o(d^{-1}).$$

Then Chebyshev's inequality gives us, for any $\tau > 0$,

$$P\left[\left|\sum_{i=1}^d \tilde{\lambda}_{i,d}z_{im}^2 - 1\right| > \tau\right] \leq \frac{E(\sum_{i=1}^d \tilde{\lambda}_{i,d}Y_i)^4}{\tau^4} \leq \frac{o(d^{-1})}{\tau^4}.$$

Summing over $d$ gives $\sum_{d=1}^\infty P[|\sum_{i=1}^d \tilde{\lambda}_{i,d}z_{im}^2 - 1| > \tau] < \infty$ and by the Borel–Cantelli lemma, we conclude that a diagonal element $\sum_{i=1}^d \tilde{\lambda}_{i,d}z_{ij}^2$ converges to 1 almost surely.

The off-diagonal elements of $nS_D$ can be expressed as $\sum_{i=1}^d \tilde{\lambda}_{i,d}z_{im}z_{il}$. Using similar arguments to those used in the diagonal case, we have

$$P\left[\left|\sum_{i=1}^d \tilde{\lambda}_{i,d}z_{im}z_{il}\right| > \tau\right] \leq \frac{E(\sum_{i=1}^d \tilde{\lambda}_{i,d}z_{im}z_{il})^4}{\tau^4} \leq \frac{o(d^{-1})}{\tau^4},$$

and again by the Borel–Cantelli lemma, the off-diagonal elements converge to 0 almost surely.

5.2. *Proofs of Lemma 1 and Theorem 2.* The proof of Theorem 2 is divided in two parts. Since eigenvectors are associated to eigenvalues, at first, we focus on asymptotic behavior of sample eigenvalues (Section 5.2.1) and then investigate consistency or strong inconsistency of sample eigenvectors (Section 5.2.2).

5.2.1. *Proof of Lemma 1.* The proof relies heavily on the following lemma. Recall that $\varphi_k(A)$ denotes the $k$th largest eigenvalue of $A$.

LEMMA 3 (Weyl's inequality). *If $A$, $B$ are $m \times m$ real symmetric matrices, then for all $k = 1, \ldots, m$,*

$$\left.\begin{array}{c}\varphi_k(A) + \varphi_m(B) \\ \varphi_{k+1}(A) + \varphi_{m-1}(B) \\ \vdots \\ \varphi_m(A) + \varphi_k(B)\end{array}\right\} \leq \varphi_k(A+B) \leq \left\{\begin{array}{c}\varphi_k(A) + \varphi_1(B), \\ \varphi_{k-1}(A) + \varphi_2(B), \\ \vdots \\ \varphi_1(A) + \varphi_k(B).\end{array}\right.$$

This inequality is discussed in Rao [17] and its use on asymptotic studies of eigenvalues of a random matrix appeared in Eaton and Tyler [6].

Since $S$ and its dual $S_D$ share nonzero eigenvalues, one of the main ideas of the proof is working with $S_D$. By our decomposition (1.1), $nS_D = Z'\Lambda Z$. We also write $Z$ and $\Lambda$ as block matrices such that

$$(5.2) \qquad Z = \begin{pmatrix} Z_1 \\ Z_2 \\ \vdots \\ Z_{p+1} \end{pmatrix}, \qquad \Lambda = \begin{pmatrix} \Lambda_1 & O & \cdots & O \\ O & \Lambda_2 & \cdots & O \\ \vdots & \vdots & \ddots & \vdots \\ O & O & \cdots & \Lambda_{p+1} \end{pmatrix},$$

where $Z_l$ is a $k_l \times n$ matrix for each $l = 1, \ldots, p+1$ and $\Lambda_l (\equiv \Lambda_{l,d})$ is a $k_l \times k_l$ diagonal matrix for each $l = 1, \ldots, p+1$ and $O$ denotes a matrix where all elements are zeros. Now, we can write

$$(5.3) \qquad nS_D = Z'\Lambda Z = \sum_{l=1}^{p+1} Z_l' \Lambda_l Z_l.$$

While $Z_l$ depends on $d = 1, \ldots, \infty$, this dependence is not explicitly shown (e.g., by subscript) for simplicity of notation.

Note that Theorem 1 implies that when the last term in equation (5.3) is divided by $d$, it converges to an identity matrix, namely,

$$(5.4) \qquad d^{-1} Z_{p+1}' \Lambda_{p+1} Z_{p+1} \xrightarrow{p} nK \cdot I_n,$$

where $K \in (0, \infty)$ is such that $(dn)^{-1} \sum_{i \in J_{p+1}} \lambda_{i,d} \to K$. Moreover, dividing by $d^{\alpha_1}$ gives us

$$nd^{-\alpha_1} S_D = d^{-\alpha_1} Z_1' \Lambda_1 Z_1 + d^{-\alpha_1} \sum_{l=2}^{p} Z_l' \Lambda_l Z_l + d^{1-\alpha_1} d^{-1} Z_{p+1}' \Lambda_{p+1} Z_{p+1}.$$

By the assumption (b), the first term on the right-hand side converges to $Z_1' C_1 Z_1$ where $C_1$ is the $k_1 \times k_1$ diagonal matrix such that $C_1 = \text{diag}\{c_j; j \in J_1\}$ and the other terms tend to a zero matrix. Thus, we get

$$nd^{-\alpha_1} S_D \Longrightarrow Z_1' C_1 Z_1 \qquad \text{as } d \to \infty.$$

Note that the nonzero eigenvalues of $Z_1' C_1 Z_1$ are the same as the nonzero eigenvalues of $C_1^{1/2} Z_1 Z_1' C_1^{1/2}$ which is a $k_1 \times k_1$ random matrix with full rank almost surely. Since $\varphi_i(A)$ is a continuous function of the entries of $A$ (see e.g., Kato [13]), we have for $i \in J_1$,

$$\varphi_i(nd^{-\alpha_1}S_D) \Longrightarrow \varphi_i(Z_1'C_1Z_1) \quad \text{as } d \to \infty$$
$$= \varphi_i(C_1^{1/2}Z_1Z_1'C_1^{1/2}).$$

Thus, we conclude that for the sample eigenvalues in the group $J_1$, $\hat{\lambda}_i/d^{\alpha_1} = \varphi_i(d^{-\alpha_1}S_D)$ converges in distribution to $\varphi_i(n^{-1}C_1^{1/2}Z_1Z_1'C_1^{1/2})$ for $i \in J_1$.

Let us focus on eigenvalues whose indices are in the group $J_2, \ldots, J_p$. Suppose we have $\hat{\lambda}_i = O_p(d^{\alpha_j})$ for all $i \in J_j$, for $j = 1, \ldots, l-1$. Pick any $i \in J_l$. We will provide upper and lower bounds on $\hat{\lambda}_i$ by Weyl's inequality (Lemma 3). Dividing both sides of (5.3) by $d^{\alpha_l}$, we get

$$nd^{-\alpha_l}S_D = d^{-\alpha_l}\sum_{j=1}^{l-1}Z_j'\Lambda_jZ_j + d^{-\alpha_l}\sum_{j=l}^{p+1}Z_j'\Lambda_jZ_j$$

and apply Weyl's inequality for the upper bound,

(5.5)
$$\varphi_i(nd^{-\alpha_l}S_D) \leq \varphi_{1+\sum_{j=1}^{l-1}k_j}\left(d^{-\alpha_l}\sum_{j=1}^{l-1}Z_j'\Lambda_jZ_j\right)$$
$$+ \varphi_{i-\sum_{j=1}^{l-1}k_j}\left(d^{-\alpha_l}\sum_{j=l}^{p+1}Z_j'\Lambda_jZ_j\right)$$
$$= \varphi_{i-\sum_{j=1}^{l-1}k_j}\left(d^{-\alpha_l}\sum_{j=l}^{p+1}Z_j'\Lambda_jZ_j\right).$$

Note that the first term vanishes since the rank of $d^{-\alpha_l}\sum_{j=1}^{l-1}Z_j'\Lambda_jZ_j$ is at most $\sum_{j=1}^{l-1}k_j$. Also note that the matrix in the upper bound (5.5) converges to a simple form

$$d^{-\alpha_l}\sum_{j=l}^{p+1}Z_j'\Lambda_jZ_j = d^{-\alpha_l}Z_l'\Lambda_lZ_l + d^{-\alpha_l}\sum_{j=l+1}^{p+1}Z_j'\Lambda_jZ_j$$
$$\Longrightarrow Z_l'C_lZ_l \quad \text{as } d \to \infty,$$

where $C_l$ is the $k_l \times k_l$ diagonal matrix such that $C_l = \text{diag}\{c_j; j \in J_l\}$.

In order to have a lower bound of $\hat{\lambda}_i$, Weyl's inequality is applied to the expression

$$d^{-\alpha_l}\sum_{j=1}^{l}Z_j'\Lambda_jZ_j + d^{-\alpha_l}\sum_{j=l+1}^{p+1}Z_j'\Lambda_jZ_j = nd^{-\alpha_l}S_D,$$

so that

$$(5.6) \quad \varphi_i\left(d^{-\alpha_l}\sum_{j=1}^{l} Z_j'\Lambda_j Z_j\right) + \varphi_n\left(d^{-\alpha_l}\sum_{j=l+1}^{p+1} Z_j'\Lambda_j Z_j\right) \leq \varphi_i(nd^{-\alpha_l}S_D).$$

It turns out that the first term of the left-hand side is not easy to manage, so we again use Weyl's inequality to get

$$(5.7) \quad \varphi_{\sum_{j=1}^{l} k_j}\left(d^{-\alpha_l}\sum_{j=1}^{l-1} Z_j'\Lambda_j Z_j\right)$$
$$\leq \varphi_i\left(d^{-\alpha_l}\sum_{j=1}^{l} Z_j'\Lambda_j Z_j\right) + \varphi_{1-i+\sum_{j=1}^{l-1} k_j}(-d^{-\alpha_l}Z_l'\Lambda_l Z_l),$$

where the left-hand side is 0 since the rank of the matrix inside is at most $\sum_{j=1}^{l-1} k_j$. Note that since $d^{-\alpha_l}Z_l'\Lambda_l Z_l$ and $d^{-\alpha_l}\Lambda_l^{1/2}Z_l Z_l'\Lambda_l^{1/2}$ share nonzero eigenvalues, we get

$$(5.8) \quad \begin{aligned} \varphi_{1-i+\sum_{j=1}^{l} k_j}(-d^{-\alpha_l}Z_l'\Lambda_l Z_l) &= \varphi_{1-i+\sum_{j=1}^{l} k_j}(-d^{-\alpha_l}\Lambda_l^{1/2}Z_l Z_l'\Lambda_l^{1/2}) \\ &= \varphi_{k_l-i+1+\sum_{j=1}^{l-1} k_j}(-d^{-\alpha_l}\Lambda_l^{1/2}Z_l Z_l'\Lambda_l^{1/2}) \\ &= -\varphi_{i-\sum_{j=1}^{l-1} k_j}(d^{-\alpha_l}\Lambda_l^{1/2}Z_l Z_l'\Lambda_l^{1/2}) \\ &= -\varphi_{i-\sum_{j=1}^{l-1} k_j}(d^{-\alpha_l}Z_l'\Lambda_l Z_l). \end{aligned}$$

Here, we use the fact that for any $m \times m$ real symmetric matrix $A$, $\varphi_i(A) = -\varphi_{m-i+1}(-A)$ for all $i = 1, \ldots, m$.

Combining (5.6)–(5.8) gives the lower bound

$$(5.9) \quad \varphi_{i-\sum_{j=1}^{l-1} k_j}(d^{-\alpha_l}Z_l'\Lambda_l Z_l) + \varphi_n\left(d^{-\alpha_l}\sum_{j=l+1}^{p+1} Z_j'\Lambda_j Z_j\right) \leq \varphi_i(nd^{-\alpha_l}S_D).$$

Note that the matrix inside of the first term of the lower bound (5.9) converges to $Z_l'C_l Z_l$ in distribution. The second term converges to 0 since the matrix inside converges to a zero matrix.

The difference between the upper and lower bounds of $\varphi_i(nd^{-\alpha_l}S_D)$ converges to 0 since

$$\varphi_{i-\sum_{j=1}^{l-1} k_j}\left(d^{-\alpha_l}\sum_{j=l}^{p+1} Z_j'\Lambda_j Z_j\right) - \varphi_{i-\sum_{j=1}^{l-1} k_j}(d^{-\alpha_l}Z_l'\Lambda_l Z_l) \to 0,$$

as $d \to \infty$. This is because $\varphi$ is a continuous function and the difference between the two matrices converges to zero matrix. Therefore, $\varphi_i(nd^{-\alpha_l}S_D)$ converges to the upper or lower bound as $d \to \infty$.

22 S. JUNG AND J. S. MARRONNow since both upper and lower bound of $\varphi_i(nd^{-\alpha_l}S_D)$ converge in distribution to same quantity, we have

$$\varphi_i(nd^{-\alpha_l}S_D) \Longrightarrow \varphi_{i-\sum_{j=1}^{l-1} k_j}(Z_l'C_lZ_l) \quad \text{as } d \to \infty.$$

(5.10)
$$= \varphi_{i-\sum_{j=1}^{l-1} k_j}(C_l^{1/2}Z_lZ_l'C_l^{1/2}).$$

Thus, by induction, we have the scaled $i$th sample eigenvalue $\hat{\lambda}_i/d^{\alpha_l}$ converges in distribution to $\varphi_{i-\sum_{j=1}^{l-1} k_j}(n^{-1}C_l^{1/2}Z_lZ_l'C_l^{1/2})$ for $i \in J_l$, $l = 1, \ldots, p$, as desired.

Now, let us focus on the rest of the sample eigenvalues $\hat{\lambda}_i$, $i = \kappa+1, \ldots, n$. For any $i$, again by Weyl's upper bound inequality, we get

$$\varphi_i(nd^{-1}S_D) \leq \varphi_{i-\kappa}(d^{-1}Z_{p+1}'\Lambda_{p+1}Z_{p+1}) + \varphi_{\kappa+1}\left(d^{-1}\sum_{j=1}^{p}Z_j'\Lambda_jZ_j\right)$$

$$= \varphi_{i-\kappa}(d^{-1}Z_{p+1}'\Lambda_{p+1}Z_{p+1}),$$

where the second term on the right-hand side vanishes since the matrix inside is of rank at most $\kappa$. Also for lower bound, we have

$$\varphi_i(nd^{-1}S_D) \geq \varphi_i(d^{-1}Z_{p+1}'\Lambda_{p+1}Z_{p+1}) + \varphi_n\left(d^{-1}\sum_{j=1}^{p}Z_j'\Lambda_jZ_j\right)$$

$$= \varphi_i(d^{-1}Z_{p+1}'\Lambda_{p+1}Z_{p+1}),$$

where the second term vanishes since $\kappa < n$. Thus, we have complete bounds for $\varphi_i(nd^{-1}S_D)$ such that

$$\varphi_i(d^{-1}Z_{p+1}'\Lambda_{p+1}Z_{p+1}) \leq \varphi_i(nd^{-1}S_D) \leq \varphi_{i-\kappa}(d^{-1}Z_{p+1}'\Lambda_{p+1}Z_{p+1})$$

for all $i = \kappa+1, \ldots, n$. However, by (5.4), the matrix in both bounds converges to $nK \cdot I_n$ in probability. Thus, lower and upper bounds of $\varphi_i(d^{-1}S_D)$ converge to $K$ in probability for $i = \kappa+1, \ldots, n$, which completes the proof.

5.2.2. *Proof of Theorem 2.* We begin by defining a standardized version of the sample covariance matrix, not to be confused with the dual $S_D$, as

$$\tilde{S} = \Lambda^{-1/2}U'SU\Lambda^{-1/2}$$

(5.11)
$$= \Lambda^{-1/2}U'(\hat{U}\hat{\Lambda}\hat{U}')U\Lambda^{-1/2}$$

$$= \Lambda^{-1/2}P\hat{\Lambda}P'\Lambda^{-1/2},$$

where $P = U'\hat{U} = \{u_i'\hat{u}_j\}_{ij} \equiv \{p_{ij}\}_{ij}$. Note that elements of $P$ are inner products between population eigenvectors and sample eigenvectors. Since $\tilde{S}$ is standardized, we have by $S = n^{-1}XX'$ and $X = U\Lambda^{1/2}Z$,

(5.12)
$$\tilde{S} = n^{-1}ZZ'.$$

Note that the angle between two directions can be formulated as an inner product of the two direction vectors. Thus, we will investigate the behavior of the inner product matrix $P$ as $d \to \infty$, by showing that

$$\sum_{j \in J_l} p_{ji}^2 \xrightarrow{p} 1 \quad \text{as } d \to \infty \tag{5.13}$$

for all $i \in J_l$, $l = 1, \ldots, p$ and

$$p_{ii}^2 \xrightarrow{p} 0 \quad \text{as } d \to \infty \tag{5.14}$$

for all $i = \kappa + 1, \ldots, n$.

Suppose for now we have the result of (5.13) and (5.14). Then for any $i \in J_l$, $l = 1, \ldots, p$,

$$\begin{aligned}
\text{Angle}(\hat{u}_i, \text{span}\{u_j : j \in J_l\}) &= \arccos\left(\frac{\hat{u}_i'[\text{Proj}_{\text{span}\{u_j : j \in J_l\}} \hat{u}_i]}{\|\hat{u}_i\|_2 \cdot \|[\text{Proj}_{\text{span}\{u_j : j \in J_l\}} \hat{u}_i]\|_2}\right) \\
&= \arccos\left(\frac{\hat{u}_i'(\sum_{j \in J_l}(u_j'\hat{u}_i)u_j)}{\|\hat{u}_i\|_2 \cdot \|\sum_{j \in J_l}(u_j'\hat{u}_i)u_j\|_2}\right) \\
&= \arccos\left(\frac{\sum_{j \in J_l}(u_j'\hat{u}_i)^2}{1 \cdot (\sum_{j \in J_l}(u_j'\hat{u}_i)^2)^{1/2}}\right) \\
&= \arccos\left(\left(\sum_{j \in J_l} p_{ji}^2\right)^{1/2}\right) \\
&\xrightarrow{p} 0 \quad \text{as } d \to \infty,
\end{aligned}$$

by (5.13) and for $i = \kappa + 1, \ldots, n$,

$$\begin{aligned}
\text{Angle}(\hat{u}_i, u_i) &= \arccos(|u_i'\hat{u}_i|) \\
&= \arccos(|p_{ii}|) \\
&\xrightarrow{p} \frac{\pi}{2} \quad \text{as } d \to \infty,
\end{aligned}$$

by (5.14), as desired.

Therefore, it is enough to show (5.13) and (5.14). We begin with taking $j$th diagonal entry of $\tilde{S}$, $\tilde{s}_{jj}$, from (5.11) and (5.12),

$$\tilde{s}_{jj} = \lambda_j^{-1} \sum_{i=1}^{n} \hat{\lambda}_i p_{ji}^2 = n^{-1} z_j z_j',$$

where $z_j$ denotes the $j$th row vector of $Z$. Since

$$\lambda_j^{-1} \hat{\lambda}_i p_{ji}^2 \leq n^{-1} z_j z_j', \tag{5.15}$$

we have at most
$$p_{ji}^2 = O_p\left(\frac{\lambda_j}{\hat{\lambda}_i}\right)$$
for all $i = 1, \ldots, n$, $j = 1, \ldots, d$. Note that by Lemma 1, we have for $i \in J_{l_1}$, $j \in J_{l_2}$ where $1 \leq l_1 < l_2 \leq p+1$,

$$(5.16) \qquad p_{ji}^2 = O_p\left(\frac{\lambda_j}{\hat{\lambda}_i}\right) = \begin{cases} O_p(d^{\alpha_{l_2} - \alpha_{l_1}}), & \text{if } l_2 \leq p, \\ O_p(d^{1-\alpha_{l_1}}), & \text{if } l_2 = p+1, \end{cases}$$

so that $p_{ji}^2 \xrightarrow{p} 0$ as $d \to \infty$ in both cases.

Note that the inner product matrix $P$ is also a unitary matrix. The norm of the $i$th column vector of $P$ must be 1 for all $d$, i.e. $\sum_{j=1}^d p_{ji}^2 = 1$. Thus, (5.13) is equivalent to $\sum_{j \in \{1,\ldots,d\}\setminus J_l} p_{ji}^2 \xrightarrow{p} 0$ as $d \to \infty$.

Now for any $i \in J_1$,
$$\sum_{j \in \{1,\ldots,d\}\setminus J_1} p_{ji}^2 = \sum_{j \in J_2 \cup \cdots \cup J_p} p_{ji}^2 + \sum_{j \in J_{p+1}} p_{ji}^2.$$

Since the first term on the right-hand side is a finite sum of quantities converging to 0, it converges to 0 almost surely as $d$ tends to infinity. By (5.15), we have an upper bound for the second term,

$$\sum_{j \in J_{p+1}} p_{ji}^2 = \sum_{j \in J_{p+1}} \lambda_j^{-1} \hat{\lambda}_i p_{ji}^2 \frac{\lambda_j}{\hat{\lambda}_i}$$
$$\leq \frac{\sum_{j \in J_{p+1}} n^{-1} z_j z_j' \lambda_j}{d} \frac{d}{\hat{\lambda}_i} = \frac{\sum_{k=1}^n \sum_{j=\kappa+1}^d z_{j,k}^2 \lambda_j}{nd} \frac{d}{\hat{\lambda}_i},$$

where the $z_{j,k}$'s are the entries of a row random vector $z_j$. Note that by applying Theorem 1 with $\Sigma_d = \text{diag}\{\lambda_{\kappa+1}, \ldots, \lambda_d\}$, we have $\sum_{j=\kappa+1}^d z_{j,k}^2 \lambda_j/d \xrightarrow{p} 1$ as $d \to \infty$. Also by Lemma 1, the upper bound converges to 0 in probability. Thus, we get

$$\sum_{j \in \{1,\ldots,d\}\setminus J_1} p_{ji}^2 \xrightarrow{p} 0 \qquad \text{as } d \to \infty,$$

which is equivalent to

$$(5.17) \qquad \sum_{j \in J_1} p_{ji}^2 \xrightarrow{p} 1 \qquad \text{as } d \to \infty.$$

Let us focus on the group $J_2, \ldots, J_p$. For any $l = 2, \ldots, p$, suppose we have $\sum_{j \in J_m} p_{ji}^2 \xrightarrow{p} 1$ as $d \to \infty$ for all $i \in J_m$, $m = 1, \ldots, l-1$. Note that it implies that for any $j \in J_m$, $m = 1, \ldots, l-1$,

$$(5.18) \qquad \sum_{i \in \{1,\ldots,d\}\setminus J_m} p_{ji}^2 \xrightarrow{p} 0 \qquad \text{as } d \to \infty,$$

since

$$\sum_{j \in J_m} \sum_{i \in \{1,\ldots,d\} \setminus J_m} p_{ji}^2 = \sum_{j \in J_m} \sum_{i=1}^{d} p_{ji}^2 - \sum_{j \in J_m} \sum_{i \in J_m} p_{ji}^2 \xrightarrow{p} \sum_{j \in J_m} 1 - \sum_{i \in J_m} 1 = 0,$$

as $d \to \infty$.

Now, pick $i \in J_l$. We have

$$\sum_{j \in \{1,\ldots,d\} \setminus J_l} p_{ji}^2 = \sum_{j \in J_1 \cup \cdots \cup J_{l-1}} p_{ji}^2 + \sum_{j \in J_{l+1} \cup \cdots \cup J_p} p_{ji}^2 + \sum_{j \in J_{p+1}} p_{ji}^2.$$

Note that the first term is bounded as

$$\sum_{j \in J_1 \cup \cdots \cup J_{l-1}} p_{ji}^2 \leq \sum_{i \in J_l} \sum_{j \in J_1 \cup \cdots \cup J_{l-1}} p_{ji}^2 \leq \sum_{m=1}^{l-1} \sum_{j \in J_m} \left( \sum_{i \in \{1,\ldots,d\} \setminus J_m} p_{ji}^2 \right) \xrightarrow{p} 0$$

by (5.18). The second term also converges to 0 by (5.16). The last term is also bounded as

$$\sum_{j \in J_{p+1}} p_{ji}^2 = \sum_{j \in J_{p+1}} \lambda_j^{-1} \hat{\lambda}_i p_{ji}^2 \frac{\lambda_j}{\hat{\lambda}_i} \leq \frac{\sum_{j \in J_{p+1}} n^{-1} z_j z_j' \lambda_j}{d} \frac{d}{\hat{\lambda}_i},$$

so that it also converges to 0 in probability. Thus, we have $\sum_{j \in \{1,\ldots,d\} \setminus J_l} p_{ji}^2 \xrightarrow{p} 0$ as $d \to \infty$ which implies that

$$\sum_{j \in J_l} p_{ji}^2 \xrightarrow{p} 1 \qquad \text{as } d \to \infty.$$

Thus, by induction, (5.13) is proved.

For $i = \kappa + 1, \ldots, n$, we have $\lambda_i^{-1} \hat{\lambda}_i p_{ii}^2 \leq n^{-1} z_i z_i'$, and so

$$p_{ii}^2 \leq \hat{\lambda}_i^{-1} \lambda_i n^{-1} z_i z_i' = O_p(\hat{\lambda}_i^{-1} \lambda_i),$$

which implies (5.14) by the assumption (c) and Lemma 1, and the proof is completed.

5.3. *Proof of Corollary 1.* The proof follows the same lines as the proof of Theorem 2, with convergence in probability replaced by almost sure convergence.

5.4. *Proof of Corollary 2.* From the proof of Theorem 2, write the inner product matrix $P$ of (5.11) as a block matrix such that

$$P = \begin{pmatrix} P_{11} & \cdots & P_{1p} & P_{1,p+1} \\ \vdots & \ddots & \vdots & \vdots \\ P_{p1} & \cdots & P_{pp} & P_{p,p+1} \\ P_{p+1,1} & \cdots & P_{p+1,p} & P_{p+1,p+1} \end{pmatrix},$$

where each $P_{ij}$ is a $k_i \times k_j$ random matrix. In the proof of Theorem 2, we have shown that $P_{ii}$, $i = 1, \ldots, p$, tends to be a unitary matrix and $P_{ij}$, $i \neq j$, tends to be a zero matrix as $d \to \infty$. Likewise, $\Lambda$ and $\hat{\Lambda}$ can be blocked similarly as $\Lambda = \text{diag}\{\Lambda_i : i = 1, \ldots, p+1\}$ and $\hat{\Lambda} = \text{diag}\{\hat{\Lambda}_i : i = 1, \ldots, p+1\}$.

Now, pick $l \in \{1, \ldots, p\}$. The $l$th block diagonal of $\tilde{S}$, $\tilde{S}_{ll}$, is expressed as $\tilde{S}_{ll} = \sum_{j=1}^{p+1} \Lambda_l^{-1/2} P_{lj} \hat{\Lambda}_l P'_{lj} \Lambda_l^{-1/2}$. Since $P_{ij} \to 0$, $i \neq j$, we get

$$\|\tilde{S}_{ll} - \Lambda_l^{-1/2} P_{ll} \hat{\Lambda}_l P'_{ll} \Lambda_l^{-1/2}\|_F \xrightarrow{p} 0$$

as $d \to \infty$, where $\|\cdot\|_F$ is the Frobenius norm of matrices defined by $\|A\|_F = (\sum_{i,j} A_{ij}^2)^{1/2}$.

Note that by (5.12), $\tilde{S}_{ll}$ can be replaced by $n^{-1} Z_l Z'_l$. We also have $d^{-\alpha_l} \Lambda_l \to C_l$ by the assumption (b) and $d^{-\alpha_l} \hat{\Lambda}_l \xrightarrow{p} \text{diag}\{\varphi(n^{-1} C_l^{1/2} Z_l Z'_l C_l^{1/2})\}$ by (5.10). Thus, we get

$$\|n^{-1} Z_l Z'_l - C_l^{-1/2} P_{ll} \text{diag}\{\varphi(n^{-1} C_l^{1/2} Z_l Z'_l C_l^{1/2})\} P'_{ll} C_l^{-1/2}\|_F \xrightarrow{p} 0$$

as $d \to \infty$.

Also note that since $n^{-1} Z_l Z'_l \to I_{k_l}$ almost surely as $n \to \infty$, we get $n^{-1} C_l^{1/2} Z_l Z'_l C_l^{1/2} \to C_l$ and $\text{diag}\{\varphi(n^{-1} C_l^{1/2} Z_l Z'_l C_l^{1/2})\} \to C_l$ almost surely as $n \to \infty$. Using the fact that the Frobenius norm is unitarily invariant and $\|AB\|_F \leq \|A\|_F \|B\|_F$ for any square matrices $A$ and $B$, we get

$\|P'_{ll} C_l P_{ll} - C_l\|_F$

$\leq \|P'_{ll} C_l P_{ll} - \text{diag}\{\varphi(n^{-1} C_l^{1/2} Z_l Z'_l C_l^{1/2})\}\|_F + o_p(1)$

$= \|C_l - P_{ll} \text{diag}\{\varphi(n^{-1} C_l^{1/2} Z_l Z'_l C_l^{1/2})\} P'_{ll}\|_F + o_p(1)$

(5.19) $\leq \|n^{-1} C_l^{1/2} Z_l Z'_l C_l^{1/2} - P_{ll} \text{diag}\{\varphi(n^{-1} C_l^{1/2} Z_l Z'_l C_l^{1/2})\} P'_{ll}\|_F + o_p(1)$

$\leq \|C_l^{1/2}\|_F^2 \|n^{-1} Z_l Z'_l$

$\quad - C_l^{-1/2} P_{ll} \text{diag}\{\varphi(n^{-1} C_l^{1/2} Z_l Z'_l C_l^{1/2})\} P'_{ll} C_l^{-1/2}\|_F + o_p(1)$

$\xrightarrow{p} 0 \quad \text{as } d, n \to \infty.$

Note that in order to have (5.19), $P_{ll}$ must converge to $\text{diag}\{\pm 1, \pm 1, \ldots, \pm 1\}$ since diagonal entries of $C_l$ are distinct and a spectral decomposition is unique up to sign changes. Let $l = 1$ for simplicity. Now for any $m = 2, \ldots, k_1$, $p_{m1}^2 \xrightarrow{p} 0$ since

$$\|P'_{11} C_1 P_{11} - C_1\|_F^2 \geq \sum_{j=1}^{k_1} (c_1 - c_j)^2 p_{j1}^2 \geq (c_1 - c_m)^2 p_{m1}^2.$$

This leads to $p_{11}^2 \xrightarrow{p} 1$ as $d, n \to \infty$. By induction, $p_{ii}^2 \xrightarrow{p} 1$ for all $i \in J_l$, $l = 1, \ldots, p$. Therefore, $\text{Angle}(\hat{u}_i, u_i) = \arccos(|p_{ii}|) \xrightarrow{p} 0$ as $d, n \to \infty.$

If the assumptions of Corollary 1 also hold, then every convergence in the proof is replaced by almost sure convergence, which completes the proof.

5.5. *Proof of Corollary 3.* With Gaussian assumption, noticing $C_l^{1/2} Z_l \times Z_l' C_l^{1/2} \sim \mathcal{W}_{k_l}(n, C_l)$ gives the first result. When $k_l = 1$, the assumption (b) and that $C_l^{1/2} Z_l Z_l' C_l^{1/2} \sim c_i \chi_n^2$ imply that

$$\frac{\hat{\lambda}_i}{\lambda_i} = \frac{\hat{\lambda}_i}{c_i d^{\alpha_l}} \cdot \frac{c_i d^{\alpha_l}}{\lambda_i} \Longrightarrow \frac{\chi_n^2}{n} \quad \text{as } d \to \infty.$$

**Acknowledgments.** The authors are very grateful to John T. Kent (University of Leeds, UK) for the insightful Example 3.1. We also thank anonymous referees for many valuable suggestions.

DEPARTMENT OF STATISTICS
AND OPERATIONS RESEARCH
UNIVERSITY OF NORTH CAROLINA
CHAPEL HILL, NORTH CAROLINA 27599
USA
E-MAIL: sungkyu@email.unc.edu
        marron@email.unc.edu